\documentclass{article}       
\usepackage{amsrefs}
\newtheorem{theorem}{Theorem}
\newtheorem{lemma}{Lemma}
\newtheorem{corollary}{Corollary}

\def\email{{\sl e-mail: }}
\date{}
\begin{document}

\title{Geodesic mappings and concircular vector fields}


\author{Igor G. Shandra}

\maketitle

\begin{abstract}
In the present paper we study geodesic mappings of special pseudo-Riemannian manifolds called  $V_n(K)$-spaces. We prove that the set of solutions of the system of equations of geodesic mappings on  $V_n(K)$-spaces $(K\neq0)$  forms a special Jordan algebra and the set of solutions generated by consircular fields is an ideal of this algebra. We show that pseudo-Riemannian manifolds admitting a concircular field of the basic type  form the class of manifolds closed with respect to the geodesic mappings.
\smallskip

\noindent{\bf Keywords:} Pseudo-Riemannian manifold, Jordan algebra, concircular fields, geodesic mappings.
\smallskip

\noindent{\bf Mathematical Subject Classification:} 53C20; 53C25; 53C40.
\footnote{I.G. Shandra, 
    Dept.\ of Data Analysis, Decision-Making and Financial Technology,
    Finan\-cial University under the Government of the Russian Federation, 
    Leningradsky Prospect 49-55, 125468 Moscow, Russia, 
      \email{ma-tematika@yandex.ru}}          

\end{abstract}

%
%
%
\def\a{\alpha}
\def\b{\beta}
\def\g{\gamma}
\def\la{\lambda}
\def\G{\Gamma}
\def\r{\varrho}
\def\an{\hbox{$A_n$}}
\def\ann{\hbox{$\bar A_n$}}
\def\mn{\hbox{$M_n$}}
\font\msbm = msbm10 
\def\N{\hbox{\msbm{N}}}
\def\vn{\hbox{$V_n$}}
\def\vnn{\hbox{$\bar V_n$}}
\def\vnk{\hbox{$V_n(K)$}}
\def\conv{\hbox{${\rm Con}(V_n)$}}
\def\na{\nabla}
\def\ph{\varphi}
\def\be#1{\begin{equation}\label{#1}}
\def\ee{\end{equation}}
\def\re#1{\eqref{#1}}
\def\eqref#1{\hbox{$(\ref{#1})$}}
\def\R{{\msbm{R}}}
\def\nad#1#2{{\buildrel{#1} \over{#2}\strut}\!}
\def\pod#1#2{{\mathrel{\mathop{#2}\limits_{#1}}\strut}\!}
%
%

\section{Introduction}
\label{intro}

The problem of geodesic mappings of pseudo-Riemannian manifold was first started by Levi-Civita [12]. There exists many monographs and papers devoted to the theory of geodesic mappings and transformations 
\cite{1,2,3,4,5,7,8,9,10,11,12,13,14,15,16,17,18,19,20,21,22,23,24,25,26,27,28,29,30,31,32,mid,mi80r,v1,v2,v3}.
Geodesic mappings play an important role in the general theory of relativity \cite{7,23}.

Let $\an=(\mn,\na)$  be a $n$-dimensional manifolds \mn\  with affine connection~$\na$  without torsion. We denote the ring of smooth functions on \mn\  by $f(\mn)$, the Lie algebra of smooth vector fields on \mn\  by  $X(\mn)$ and arbitrary smooth vector fields on \mn\  by  $X,Y,Z$. 

A diffeomorphism  $f\colon\an\to\ann$ is called a {\it geodesic mapping} of \an\   onto \ann\  if $f$  maps any geodesic curve on \an\  onto a geodesic curve on \ann\  \cite{21,30}.

A manifold \an\  admits a geodesic mapping onto  \ann\ if and only if the equation~\cite{21,30}
$$
\bar\na_XY=\na_XY+\psi(X)Y+\psi(Y)X
$$
holds for any vector fields  $X,Y$   and where $\psi$  is a differential form on  $\mn (=\bar M_n)$. 

If $\psi$  then geodesic mapping is called {\it trivial} and {\it nontrivial} if  $\psi\neq0$.

Let $\vn=(\mn,g)$  be an $n$-dimensional pseudo-Riemannian manifolds with metric tensor $g$  and $\na$  be a {\it Levi-Civita connection}.

A pseudo-Riemannian manifold \vn\  admits a geodesic mapping onto pseudo-Riemannian manifold  \vnn\ if and only if there exists a differential form   on \vn\  such that the {\it Levi-Civita equation} \cite{21,30}
\be{1}
(\na_Z\bar g)(X,Y) = 2\psi(Z)\bar g(X,Y)+\psi(X)\bar g(Y,Z)+\psi(Y)\bar g(X,Z)
\ee
holds for any vector fields  $X,Y,Z$. Or in the coordinate form
\be{2}
\bar g_{ij,k}=2\psi_k\bar g_{ij}+\psi_i\bar g_{jk}+\psi_j\bar g_{ik},
\ee
where  $\psi_i=\na_i\Psi$, $\Psi$  is a scalar field.
The Levi-Civita equations \re{1} is not linear so that is not convenient for investigations. Sinyukov \cite{21,30} proved that a pseudo-Riemannian manifold  \vn\ admits a geodesic mapping if and only if there exist a differential form $\la(X)$  and a regular symmetric bilinear form $a(X,Y)$  
on~\vn\ such that the equation
\be{3}
(\na_Z a)(X,Y) = \la(X) g(Y,Z)+\la(Y) g(X,Z)
\ee
holds for any vector fields  $X,Y,Z$. Or in the coordinate form
\be{4}
a_{ij,k}=\la_ig_{jk}+\la_j g_{ik},
\ee
where   $\la_i=-a^s_i\psi_s$, $a^s_i=g^{st}a_{ti}$, $g^{st}$   are the contravariant components of the metric  $g$. 
Note that  $\la_i=\na_i\Lambda$,    $\Lambda$ is a scalar field.

If \vn\ $(n>2)$  admits two linearly independent solutions not proportional to the metric tensor  $g$ then \cite{21}
\be{5}
(\na_Y\la)(X)=K\,a(X,Y)+\mu\,g(X,Y),
\ee
\\[-10mm]
\be{6}
\na_X\mu=2K\,\la(X),
\ee
where $K$  is a constant and  $\mu$ is a scalar field on  \vn. Or in the coordinate form
\be{7}
\na_j\la_i=K\,a_{ij}+\mu\,g_{ij},
\ee
\\[-10mm]
\be{8}
\na_k\mu=2K\,\la_k.
\ee
A pseudo-Riemannian manifold satisfying the equations \re{3}, \re{5}, \re{6} is called a  {\it \vnk-space}.

This spaces for Riemannian manifolds were introduced by Solodovnikov \cite{31} 
 as $V(K)$-space  and with another problem for pseudo-Riemannian manifolds were introduced by Mikes \cite{mid,21} as $V_n(B)$-space (in this case $B=-K$).

A vector field  $\ph$ on a pseudo-Riemannian manifold \vn\ is called a {\it concircular}~if  
\be{9}
(\na_Y\ph)X=\r\,g(X,Y),
\ee
where  $\r$ is a scalar field on  \vn, see~Yano \cite{yaco}.

If $\r\neq0$  a concircular field belongs to the {\it basic type} and belongs to the {\it exceptional type} otherwise.

A pseudo-Riemannian manifold \vn\  admitting a concircular field is called an {\it equidistant space} \cite{21,30}. 
The equidistant space belongs to the basic type if it admits a concircular field of the mane type and belongs to the exceptional type if it admits concircular fields only the exceptional type.

Concircular fields play an important role in the theories of conformal and geodesic mappings and  transformations. They were studied by a number of geometers: Brinkmann \cite{bri}, Fialkow \cite{6}, Yano \cite{yaco}, Sinyukov \cite{30}, Aminova \cite{2}, 
Mike\v s \cite{mi80r,13},  Shandra \cite{25,26,27,28}, etc. 

The linear space of all concircular fields on \vn\  denotes by  \conv.
If $\nad{1}\ph,\dots,\nad{m}\ph$  is a basis in \conv\  then the tensor field 
$$
a=\sum_{\a,\b=1}^m \pod{\a\b}C\ (\nad{\a}\ph\otimes\nad{\b}\ph)
$$
is a solution of the system \re{2}, where  \smash{$\pod{\a\b}C\ (=\pod{\b\a}C)$} are some constants. So \vn\  admits the geodesic mapping.

Pseudo-Riemannian manifolds admitting concircular fields form the class of manifolds is closed with respect to the geodesic mappings \cite{21,30}. Let pseudo-Riemannian manifold \vn\  admits a geodesic mapping onto pseudo-Riemannian manifold \vnn\  if there exists a concircular field $\ph$  on  \vn\ then there exists a concircular field  $\bar\ph$  on  \vnn\    such that
\be{10}
\bar\r = \exp(\Psi)\ (\r+g^{ij}\ph_i\psi_j).
\ee
A concircular field  $\ph$ is said to be {\it special} if \cite{21,28}
\be{11}
Z(\r) = K\,g(Z,\ph),
\ee
where $\ph$  is a constant, and is said to be is said to be {\it convergent} \cite{29} if  $\r$ is a constant.
A pseudo-Riemannian manifold \vn\  admitting a convergent field is called a {\it Shirokov space}.

If there exist two linearly independent concircular field on \vn\  then all concircular fields on \vn\  are special with the same constant $K$, see \cite{21}. A pseudo-Riemannian manifold \vn\  admitting a special concircular field is a  \vnk-space. On a  \vnk-space any concircular field is special.

\section{Shirokov spaces and \vnk\  spaces  $(K\neq0)$} 

\begin{lemma}
Let pseudo-Riemannian manifold $V_{n+1}=(M_{n+1},G)$  admits a convergent fields  $\tilde \ph$ such that
\be{12}
\hbox{a) \ \ } \|\tilde \ph\|<0 \hbox{\ \ and \ \ \ b) \ \ }
(\tilde\nabla_{\tilde Y}\tilde \ph )\tilde X=K\,G(\tilde X,\tilde Y),
\ee
for any vector fields  $\tilde X,\tilde Y$ on  $M_{n+1}$,  where $K\ (\neq0)$  is a constant. Then there exists the adapted coordinate system $(x^I)=(x^0,x^i)$  in which the components  $G_{IJ}$ of the metric  $G$ reduce to the form
\be{13}
G_{IJ}=\exp(2\,K\,x^0)\ \left(
\begin{array}{cc} -1 & 0 \\[2mm]
0 & \displaystyle\frac{g_{ij}(x^k)}{K}
\end{array}
\right)
\ee
where  $g_{ij}(x^k)$ is the components of the metric of some  $\vn=(\mn,g)$,  $I,J,\dots=1,\dots,n+1$, $i,j,\dots=1,\dots,n$. 
\end{lemma}

\noindent {\it Proof.} Let  $\tilde {\ph}{}^I$ be the components  of the vector fields  $\tilde {\ph}$  $g$-conjugate with a convergent fields  $\tilde {\ph}$ in a coordinate system $(x^I)$  on  $V_{n+1}=(M_{n+1},G)$. Then due to (\ref{12}b) they satisfy
\be{14}
\tilde \na_J\tilde \ph^I=K\,\delta^I_J .
\ee
Let  $D$ be the linear space of all vector fields on $V_{n+1}$  which are orthogonal to $\nad*\ph$. It easy to check that $D$  is involutive. So if we use as a natural basis of $X(M_{n+1})$  the basis  $\{e_I\}=\{\nad*\ph,e_i\}$, where $\{e_i\}$,  is the basis in  $D$, we get the coordinate system $(x^I)=(x^0,x^i)$  in which
\be{15}
\hbox{a) \ \ }  \tilde \ph{}^I=\delta^I_0;  \hbox{\ \ \ b) \ \ }   
G_{i0}=0.  
\ee 
In these coordinates the equations \re{14} are equivalent to
\be{16}
\tilde \G^I_{0J}=K\ \delta^I_J ,    
\ee
where $\tilde \G{}^I_{JK}$  are the components of the Levi-Civita connection of the metric  $G$.

Let us consider the conditions \re{16}.  If  $I=0, J=j$  we have
\be{17}
\partial_jG_{00}=0.
\ee
If  $I=0, J=0$  we get
\be{18}
\partial_jG_{00}=2K\,G_{00}.
\ee
It follows from \re{17} and \re{18} that
\ $G_{00}=C\cdot\exp{(2K\,x^0)}$,
where $C$  is a constant. Due to (\ref{12}a) that  $C<0$. We can take it such that  $C=-1$. So
\be{19}
G_{00}=- \exp{(2K\,x^0)}.
\ee
If  $I=i,J=j$ we obtain \
$
\partial_0G_{ij}=2K\,G_{ij}
$. \
So
\be{20}
G_{ij}=\exp{(2K\,x^0)}\ \frac{g_{ij}(x^k)}{K}.
\ee
It follows  from (\ref{15}b), \re{19}, \re{20} that in the coordinate system   $(x^I)=(x^0,x^i)$ components $G_{IJ}$  reduces to the form \re{13}.

Conversely, if the components $G_{IJ}$  of the metric $G$  in the coordinate system  $(x^I)=(x^0,x^i)$ reduce to the form \re{13} then the components $\tilde \G^I_{JK}$  of the Levi-Civita connection reduce to the form:
\be{21}
\tilde \G^0_{00}=K, \ \
\tilde \G^0_{0j}=0, \ \
\tilde \G^i_{0j}=\delta^i_j, \ \
\tilde \G^0_{ij}=g_{ij}, \ \
\tilde \G^k_{ij}=\G^k_{ij},
\ee
where  $\G^k_{ij}$ are the components of the Levi-Civita connection of the metric  $g$. Using direct calculations it easy to verify that vector field with components 
$\tilde \ph{}^I_0=\delta^I_0$
  by virtue \re{21} satisfies the conditions (\ref{12}a), \re{14}.\smallskip 
	
\noindent {\bf Remark 1} \ The components  $G^{IJ}$ of the inverse metric $G$  in the adapted coordinate system $(x^I)=(x^0,x^i)$  reduce to the form
\be{22}
G^{IJ}=\exp(-2\,K\,x^0)\ \left(
\begin{array}{cc} -1 & 0 \\[2mm]
0 & K\,g^{ij}(x^k)
\end{array}
\right)
\ee
%
\begin{lemma}\label{le2}
The pseudo-Riemannian manifold $V_{n+1}=(M_{n+1},G)$  with the metric defined by the conditions \re{13} admits an absolutely parallel convector field $\tilde \ph$  if and only if its components in the adapted coordinate system  $(x^I)=(x^0,x^i)$ reduce to the form
\be{23}
\tilde \ph_ I=\exp(Kx^0)\left(\r(x^k),\ph_i(x^k)\right),
\ee
where $\r(x^k)$  and  $\ph_i(x^k)$ satisfy the following equations on  $\vn=(\mn,g)$:
 \be{24}
\na_j\ph_i=\r\,g_{ij},
\ee
\be{25}
\na_j\r=K\,\ph_j.
\ee 
\end{lemma}
{\it Proof.} Let $\tilde \ph_I$   be the components of an absolutely parallel covector field  $\tilde \ph$ in the adapted coordinate system 
$(x^I)=(x^0,x^i)$  on  $V_{n+1}=(M_{n+1},G)$. So
\be{26}
\tilde \na_J\tilde \ph_I=0
\ee
If $I=0, J=0$  we get from \re{26} by virtue \re{21}
$$
\partial_0\tilde \ph_0-K\ \tilde \ph_0 =0.
$$
 
Thus
\be{27}
\tilde \ph_0=\exp(Kx^0)\,\r(x^k).
\ee
If  $I=i, J=0$:
$$
\partial_0\tilde \ph_i-K\ \tilde \ph_i =0.
$$
Hence,
\be{28}
\tilde \ph_i=\exp(Kx^0)\,\tilde \ph_i(x^k).
\ee
If  $I=0, J=j$:
 $$
\partial_j\tilde \ph_0-K\ \tilde \ph_j =0.
$$
Due to \re{27}, \re{28} we have 
$$
\na_j\r=K\,\ph_j
$$

If  $I=i, J=j$:
 $$
\partial_j\tilde \ph_i-g_{ij}\tilde \ph_0-\G^a_{ij}\tilde \ph_a=0.
$$
Thus,
 $$
\na_j\ph_i=\r\,g_{ij}.
$$

Conversely, using direct calculations it easy to check that if the covector field $\tilde \ph$  has components  
$
\tilde \ph_i=\exp(Kx^0)\ (\r(x^k),\ph_i(x^k))
$
 in the adapted coordinate system $(x^I)=(x^0,x^i)$  on 
$V_{n+1}=(M_{n+1},G)$ with metric \re{13}, where   and   satisfy the equations \re{24}, \re{25} on  , then due to \re{21}  is absolutely parallel. \smallskip

\noindent{\bf Remark 2} \ The equations \re{24}, \re{25} are the coordinate forms of the equations \re{9}, \re{11} defining a special concircular field. So the conditions \re{23} establish a one-to-one correspondence between absolutely parallel covector fields on the Shirokov space $V_{n+1}=(M_{n+1},G)$  and special concircular fields on the  \vnk-space  $K\neq0$.

In a similar way, it is possible to prove the following statement. 
\begin{lemma}\label{le3}
The pseudo-Riemannian manifold $V_{n+1}=(M_{n+1},G)$  with  the metric defined by the conditions \re{13} admits an absolutely parallel symmetric bilinear form  $\tilde a$ if and only if its components in the adapted coordinate system $(x^I)=(x^0,x^i)$  reduce to the form
\be{29}
\tilde a_{IJ}=\exp(2\,K\,x^0)\ \left(
\begin{array}{cc} \mu(x^k) & \la_i(x^k) \\[2mm]
\la_j(x^k) & a_{ij}(x^k)
\end{array}
\right)
\ee
where  $a_{ij}(x^k)$,  $\la_i(x^k)$, and  $\mu(x^k)$ satisfy the equations \re{4}, \re{7}, \re{8} on $\vn=(\mn,g)$.
\end{lemma}
{\bf Remark 3} \  The equations \re{4}, \re{7}, \re{8} define a  \vnk-space. So the conditions \re{29} establish a one-to-one correspondence between absolutely parallel symmetric bilinear forms on the Shirokov space  $V_{n+1}=(M_{n+1},G)$ and solutions of the system \re{4}, \re{7}, \re{8} defining geodesic mappings of the  \vnk-space  $(K\neq0)$.\smallskip

\noindent{\bf Remark 4} \  The set of absolutely parallel symmetric bilinear forms  on $\vn=(\mn,g)$  is special Jordan algebra $J_0$  with the operation of multiplication  
$\nad1A*\nad2A=\{\nad1A;\nad2A\}$, where $A$  is the linear operator  $g$-conjugate with a bilinear form  $a$, defined by  
$g(AX,Y)=a(X,Y)$, and  $\{\nad1A;\nad2A\}$ is a Jordan brackets 
\be{30}
\{\nad1A;\nad2A\}=\frac12\ \left( \nad1A\,\nad2A +\nad2A\,\nad1A \right).
\ee
The condition \re{30} can be rewritten in the vector form as
\be{31}
2\,\{\nad1a;\nad2a\}(X,Y)= \nad1a\left( \nad2A X,Y \right)+\nad1a\left( \nad2A Y,X \right).
\ee
Or in the coordinate form
 \be{32}
2\,\{\nad1a;\nad2a\}_{ij}= g^{ab}\ \nad1a\left(\nad1a_{ai}\nad2a_{bj} +\nad1a_{aj}\nad2a_{bi} \right).
\ee
This statement follows from the Lemma \ref{le2}.
\begin{theorem}\label{th1}
The set of solutions of the system \re{4}, \re{7}, \re{8} on a  
\vnk-space $(K\neq0)$  forms a special Jordan algebra  $J$ with the operation of multiplication  
$\left\{
(\nad1a,\nad1\la,\nad1\mu);(\nad2a,\nad2\la,\nad2\mu)=
(\nad3a,\nad3\la,\nad3\mu)
\right\}$, where
\be{33}
2\,\nad3a(X,Y)= K\left(\nad1a ( \nad2A X,Y )+\nad1a( \nad2A Y,X )\right)-
\left( \nad1\la\otimes\nad2\la + \nad2\la\otimes\nad1\la \right) (X,Y),
\ee
\be{34}
2\,\nad3\la(X)= K\left(\nad1\la (\nad2A X)+\nad2\la( \nad1A X)\right)-
\left( \nad1\mu\nad2\la(X) + \nad2\mu\nad1\la(X) \right) ,
\ee
\be{35}
\nad3\mu= K\,g^{-1}\left(\nad1\la \,\nad2\la\right)-
\nad1\mu\,\nad2\mu.
\ee
The algebra $J$  is isomorphic to the special Jordan algebra $J_0$  of absolutely parallel symmetric bilinear forms on the Shirokov space $V_{n+1}=(M_{n+1},G)$  with the metric \re{13}. 
\end{theorem}
Proof of the theorem immediately follows from the Lemma \ref{le2} and \re{22}, \re{29},~\re{32}.\\[2mm]
{\bf Remark 5} \ Due to \re{31} the unit of the algebra $J_0$  is $G$  so the unit of the algebra~$J$  is  
$\displaystyle\left(\frac gK,0,-1\right)$.
\\
{\bf Remark 6} \  If there exists a convergent fields  $\tilde \ph$ on 
 $V_{n+1}=(M_{n+1},G)$ such that  $\|\tilde \ph\|>0$. Then there exists the adapted coordinate system $(x^I)=(x^0,x^i)$  in which the components $G_{IJ}$  of the metric  $G$ reduce to the form
$$
G_{IJ}=\exp(2\,K\,x^0)\ \left(
\begin{array}{cc} 1 & 0 \\[2mm]
0 & \displaystyle\frac{-g_{ij}(x^k)}{K}
\end{array}
\right),
$$
where  $g_{ij}(x^k)$ is the components of the metric of some  
$\vn=(\mn,g)$. Using this metric and \re{31} we can define new operation of multiplication  $\{\cdot,\cdot\}_2$. It is obvious that  $\{\nad1A;\nad2A\}=-\{\nad1A;\nad2A\}_2$.
\begin{corollary}
Let  $\vn=(\mn,g)$ be a  \vnk-space $(K\neq0)$  then there exists the solution  $(a,\la,\mu)$ of the system \re{4}, \re{7}, \re{8} satisfying the following conditions: 
\be{36}
K\,a(AX,Y)-(\la\otimes\la)(X,Y)=\smash{\frac{e\,g(X,Y)}{K}},
\ee
 \be{37}
K\,\la(AX)-\mu\,\la(X)=0,
\ee
 \be{38}
K\,g^{-1}(\la,\la)-\mu^2=-e,
\ee
where $e$  takes values  $\pm1,0$.
\end{corollary}
{\it Proof.} Let $\tilde b$  be an absolutely parallel symmetric bilinear form on the Shirokov space  $V_{n+1}=(M_{n+1},G)$  with the metric \re{13}. Then as it has shown in \cite{10} there exists the absolutely parallel symmetric bilinear form  $\tilde a$ on  
$V_{n+1}=(M_{n+1},G)$  such that  $\tilde A{}^2=e$ or in the equivalent form
\be{39}
\tilde a(\tilde A\tilde X,\tilde Y)=e\,G(\tilde X,\tilde Y).
\ee
The equation \re{39} means that  $\{\tilde a,\tilde a\}=e\,G$. Hence if $(a,\la,\mu)$  is the corresponding solution of the system \re{4}, \re{7}, \re{8} on the  \vnk-space $(K\neq0)$  then taking into account \re{33}, \re{34}, \re{35} we get \re{36}, \re{37}, \re{38}. \smallskip

As mentioned above concircular fields generate a solution of the equation~\re{2}. Denote this set of solutions by  $J_c$.
\begin{theorem}\label{th2} 
$J_c$ is an ideal of  $J$.
\end{theorem}
{\it Proof.} To prove that $J_c$  is an ideal of $J$  on  $\vn=(\mn,g)$ is equivalent to prove that  $J_{0c}$ is an ideal of  $J_0$ on  $V_{n+1}=(M_{n+1},G)$, where $J_{0c}$  is the set of absolutely parallel symmetric bilinear forms generated by absolutely parallel convector fields.

Let $\nad1\ph,\dots,\nad m\ph$  be a basis of the linear space  
Conv($V_{n+1}$) of absolutely parallel convector fields on  $V_{n+1}=(M_{n+1},G)$. Then any absolutely parallel symmetric bilinear forms generated by absolutely parallel convector fields has the components
$$
\tilde b_{IJ}=\sum_{\a,\b=1}^m \pod{\a\b}C (\nad\a\ph_I\nad\b\ph_J),
$$
where $\pod{\a\b}C\ (=\pod{\b\a}C)$ are some constants. 
Let 
$\tilde a_{IJ}$  be the components of arbitrary absolutely parallel symmetric bilinear form  $\tilde a$. We should prove that  
\hbox{$\{\tilde a,\tilde b\}\in J_{0c}$}. We have 
\be{40}
2\{\tilde a,\tilde b\}=\!G^{DT}\!\!\!\sum_{\a,\b=1}^m \!\!\pod{\a\b}C (\nad\a\ph_I\nad\b\ph_D\,\tilde a_{TJ}+\nad\a\ph_J\nad\b\ph_D\,\tilde a_{TI})=\!\!\!\!
\sum_{\a,\b=1}^m \!\!\pod{\a\b}C (\nad\a\ph_I\nad\b\Phi_J+\nad\a\ph_J\nad\b\Phi_I),
\ee
where $\nad\b\Phi_I=\nad\b\ph_D\, \tilde a_{TI} \,G^{DT}$  is an absolutely parallel convector field. Therefore, 
\be{41}
\nad\b\Phi_{I}=\sum_{\g=1}^m F^\b_\g\ \nad\g\ph_I
\ee
where $F^\b_\g$  are some constants. It follows from \re{40}, \re{41} that
 $$
2\,\{\tilde a,\tilde b\}_{IJ}=\sum_{\a,\b,\g=1}^m 
\left(F^\g_\b \pod{\a\g}C + F^\g_\a \pod{\b\g}C\right)\, \nad\a\ph_I\nad\b\ph_J.
$$
Thus,  $\{\tilde a,\tilde b\}\in J_{0c}$.

\section{$V_n(0)$-spaces}

Let $(\mn.g)$  be a  $V_n(0)$-space then there exists a solution of the system
\be{42}
\na_k a_{ij}=\la_i g_{jk} + \la_j g_{ik},
\ee
\be{43}
\na_k \la_i=\mu\, g_{ik},
\ee
where $\mu$ is a constant, and $\la_i=\na_i\Lambda$. Thus, a  $V_n(0)$-space is a Shirokov space.
%
\begin{lemma}\label{le4}
If the  $V_n(0)$-space does not admit any convergent fields of the basic type and $\ph$  is an absolutely parallel convector field on it. Then there exists the sequence of absolutely parallel covector fields $\left\{\nad\a\ph\right\}\ (\a\in \N)$  such that
\be{44}
\hbox{a) \ \ }
\nad{\a+1}\ph(X)=\nad{\a}\ph(AX)-\nad\a f\,\la(X),
\hbox{\ \ \ \ b) \ \ } 
\nad\a\ph(\la^*)=0, \ \ \forall\a\in \N,
\ee                 
where  $\nad1\ph=\ph$,  $d\nad\a f=\nad\a\ph$, $\la^*$  is the vector field  $g$-conjugate with  $\la$.
\end{lemma}
{\it Proof.} Taking into account that the $V_n(0)$  does not admit any convergent fields of the basic type we obtain from \re{43} that
\be{45}
\na_k\la_i=0.
\ee

Let $\ph_i$  be the components of an absolutely parallel convector field $\ph$  on a  $V_n(0)$. Denote  $\nad1\ph=\ph$. Consider the covector field
\be{46}
\nad2\ph_i=a^t_i\,\nad1\ph-\nad1f\la_i
\ee
where  $a^t_i$ are  components of the linear operator    $A \ (a^j_i=g^{jl}a_{il})$. It follows from \re{46} due to \re{42}, \re{45}
\be{47}
\na_k\nad2\ph_i = \nad1\ph_t\la^t g_{ik},
\ee
where  $\la^t=g^{ti}\la_i$. According to our assumption it follows from \re{47} that
$$
\nad1\ph_t\la^t=0
\hbox{ \ \ and \ \ }
\na_k\nad2\ph_i=0.
$$
Applying now similar argumentation to the covector  $\nad2\ph_i$ and continuing the process in this way, we obtain the desired sequence.\smallskip

\noindent{\bf Remark 7} \ The equation (\ref{44}b) due to (\ref{44}a) can be rewritten as
\be{48}
\ph(\nad{\a-1}\la^*)=0, \ \ \forall\a\in \N,
\ee
where $\nad\a A$  is the  $\a$-s power of linear operator  $A$. 
%
\begin{theorem}\label{th3}
Let a pseudo-Riemannian manifold \vn\  be a  $V_n(0)$-space. Then there exists a convergent field of the basic type on \vn\  or there exists the sequence of linearly independent absolutely parallel convector fields $\{\nad\a\la\}$, $(\a=1,2\dots,p\leq n-1)$  such that
\be{49}
a) \ \  \nad{\a+1}\la(X)=\nad\a\la(AX)-\nad\a\Lambda\,\a(X), \ \ \ 
b) \ \  \la(\nad{\a-1}A\,\la^*)=0,\ \ \forall\a\in A,
\ee
\be{50}
\nad p\la(AX)=\nad p\Lambda\,\la(X),
\ee
where  $\nad1\la=\la$, $\la^*$    is the vector field  $g$-conjugate with  $\la$.
\end{theorem}
\noindent{\it Proof.} \ 1) It follows from \re{43} that if  $\mu\neq0$ then $\la$  is a convergent field of the basic type on  $V_n(0)$.

2) Let $\mu=0$  then  $\na\la=0$. According to the Lemma \ref{le4} and the Remark 7 we can construct the sequence of absolutely parallel convector fields $\{\nad\a\la\}\ (\a\in \N)$  such that
$$
a) \ \ \nad{\a+1}\la(X)=\nad\a\la(AX)-\nad\a\Lambda\ \la(X), \ \ \ \ 
b)  \ \ \la(\nad{\a-1}A\la^*)=0, \ \ \forall\a\in\N.
$$
This sequence contains no more than  $p\ (\leq n-1)$ linearly independent covectors. Otherwise, $V_n(0)$  will be locally flat and so it will admit a convergent field of the basic type. Thus, 
$$
\nad{p+1}\la=\sum_{\a=1}^p C_\a\nad\a\la,
$$
where $C_\a$  are constants and $\nad1\la,\dots,\nad p\la$  are linearly independent. Changing $\nad\a\Lambda$  (defined to a constant) we can make  $\nad{p+1}\la=0$. So we get \re{50}. 
%
\begin{corollary}\label{co2}
 If the  $V_n(0)$-space does not admit any converging fields of the basic type and  $\ph$ is an absolutely parallel convector field on it. Then 
\be{51}
\nad{\a-1}\la(\ph^*)=0, \ \ \forall\a\in\N
\ee
where $\ph^*$  is the vector field  $g$-conjugate with  $\ph$.
\end{corollary}
{\it Proof.} \ We get from \re{48}
$$
(\nad{\a-1}A\la^*)=\nad{\a-1}A\la(\ph^*)=\nad{\a-1}\la(\ph^*)=0.
$$
The following statement holds.
%
\begin{theorem}\label{th4}
Let pseudo-Riemannian manifold \vn\  admits a geodesic mapping onto pseudo-Riemannian manifold \vnn\  if there exists a concircular field of the basic type on \vnn\  then there exists a concircular field of the basic type on  \vn.
\end{theorem}
{\it Proof.} \ Let $\bar\ph$  be a concircular field of the basic type on  \vnn\ ($\bar\r\neq0$)   then there exists a concircular field $\ph$ on  \vn. Let us suppose the contrary that  \vn\  does not admit concircular fields of the basic type. It means that  $\r=0$. So $\ph$  is an absolutely parallel convector field and, therefore, \vn\  is a  $V_n(0)$-space \cite{27}. So according to Theorem \ref{th3} there exists \vn\  on the sequence of linearly independent absolutely parallel convector fields 
$\{\nad\a\la\}\ (\a=1,2,\dots,p\leq n-1)$
 satisfying \re{49}, \re{50}. The equation \re{50} in the coordinate form can be written as
\be{52}
a^t_i\nad p\la_t=\nad p\Lambda\la_i.
\ee
Contracting \re{52} with  $\bar a{}^i_j$ (the inverse operator to  $a^i_j$) by $i$  and taking into account that 
$\la_i=-a^t_i\psi_t$
  we get
\be{53}
\nad p\la_j=-\nad p\Lambda\psi_j .
\ee
The condition \re{51} means that  $\ph^t\nad p\la_t=0$. Hence, due to  $\nad p\Lambda\neq0$ it follows from \re{53} that
$
\ph^t\psi_t=0.
$
From another hand since $\bar\r\neq0$  and  $\r=0$ the equation 
\re{10} gives us
$
\ph^t\psi_t\neq 0.
$
This contradiction proves the theorem. 
\\[1mm]
{\bf Remark 8} \  The Theorem \ref{th4} shows that pseudo-Riemannian manifolds admitting a concircular field of the basic type (i.e. equidistant spaces of the basic type) form the class of manifolds closed with respect to the geodesic mappings. The same properties have spaces of constant curvature \cite{21,30}, Einstein spaces \cite{14,21}, and  \vnk-spaces \cite{21}.
%
\begin{corollary}\label{co3}
Let an equidistant space of the exeptional type \vn\  admits a geodesic mapping onto a pseudo-Riemannian manifold \vnn\  then \vnn\  is a equidistant space of the exeptional type.
\end{corollary}


\end{document}